\documentclass{article}\begin{document}
\centerline{\bf\Large Generalized Elliptic Integrals and Applications}
\[
\]
\centerline{\bf N.D. Bagis}
\centerline{email: nikosbagis@hotmail.gr}
\[
\]
 
\[
\]
\centerline{\bf Abstract}

\begin{quote}
We use some general properties, presented in previous work, to evaluate special cases of integrals relating Rogers-Ramanujan continued fraction, eta function and elliptic integrals.  

\end{quote}

\section{Introductory definitions and formulas}

Let 
\begin{equation}
\eta(\tau):=e^{\pi i\tau/12} \prod^{\infty}_{n=1}(1-e^{2\pi in\tau})
\end{equation}
denotes the Dedekind eta function which is defined in the upper half complex plane. It not defined for real $\tau$.
\begin{equation}
{}_2F_{1}[a,b,c;x]:=\sum^{\infty}_{n=0}\frac{(a)_n (b)_n}{(c)_n}\frac{x^n}{n!}
\end{equation}
is the Gauss Hypergeometric function, $ (s)_k:=\frac{\Gamma(s+k)}{\Gamma(s)}$ and $\Gamma(s):=\int^{\infty}_{0}e^{-t}t^{s-1}dt$ is Euler's Gamma Function.
\\
For $\left|q\right|<1$, the Rogers Ramanujan continued fraction (RRCF) is defined as
\begin{equation}
R(q):=\frac{q^{1/5}}{1+}\frac{q^1}{1+}\frac{q^2}{1+}\frac{q^3}{1+}...  
\end{equation}
We also for $|q|<1$ define
\begin{equation}
f(-q):=\prod^{\infty}_{n=1}(1-q^n)
\end{equation}
to be the eta function of Ramanujan. Also hold the following relations proved by Ramanujan
\begin{equation}
\frac{1}{R(q)}-1-R(q)=\frac{f(-q^{1/5})}{q^{1/5}f(-q^5)} 
\end{equation} 
\begin{equation}
\frac{1}{R^5(q)}-11-R^5(q)=\frac{f(-q)^6}{q f(-q^5)^6} 
\end{equation}
For the derivative of the (RRCF) (see [5]) holds 
\begin{equation}
R'(q)=5^{-1}q^{-5/6}f(-q)^4R(q)\sqrt[6]{R(q)^{-5}-11-R(q)^5}
\end{equation} 
Some interesting results can be found if we consider the Appell $F_1$ hypergeometric function. This function is defined as 
\begin{equation}
F_1[a,b_1,b_2,c;x,y]:=\sum^{\infty}_{m=0}\sum^{\infty}_{n=0}\frac{(a)_{m+n}}{(c)_{m+n}}\frac{(b_1)_m (b_2)_{n}}{m! n!}x^m y^n
\end{equation}
It is known that
\begin{equation}
\int x^m(ax^2+bx+c)^n dx=\frac{c^n x^{m+1}}{m+1}F_1\left[m+1,-n,-n,m+2;\frac{x}{\rho_1},\frac{x}{\rho_2}\right] 
\end{equation}
where $\rho_{1,2}=\frac{b\pm \sqrt{b^2-4ac}}{2a}$.

\section{Evaluations of elliptic integrals and Rogers-Ramanujan continued fraction}
 
From relation (7) we have multiplying both sides with $G(R(q))$ and integrating 
\begin{equation}
\int^{b}_{a}f(-q)^4q^{-5/6}G(R(q))dq=5\int^{R(b)}_{R(a)}\frac{G(x)}{x\sqrt[6]{x^{-5}-11-x^5}}dx
\end{equation}
Setting $a=0$,  $b=R^{-1}\left(\sqrt[5]{\frac{-11+5\sqrt{5}}{2}}\right)$ in (10) and  making the change of variables $x\rightarrow w^{1/5}$, $w\rightarrow y/\rho_2$ with $\rho_2=\frac{-11+5\sqrt{5}}{2}$, and $w^{-1}-11-w\rightarrow t$ and then using the formula (see [13])
\begin{equation}
{}_2F_1\left(a,b,c;z\right)=\frac{\Gamma(c)}{\Gamma(b)\Gamma(c-b)}\int^{1}_{0}t^{b-1}(1-t)^{c-b-1}(1-zt)^{-a}dt
\end{equation}
\\ 
we get the following theorem\\
\\
\textbf{Theorem 1.}\\
For every $\nu>0$
$$
\int^{R^{(-1)}\left(\sqrt[5]{\frac{-11+5\sqrt{5}}{2}}\right)}_{0}f(-q)^4q^{-5/6}R(q)^{5\nu}dq=
$$
\begin{equation}
=\Gamma\left(\frac{5}{6}\right)\left(\frac{11+5 \sqrt{5}}{2}\right)^{-\frac{1}{6}-\nu} \frac{\Gamma\left(\frac{1}{6}+\nu\right)}{\Gamma(1+\nu)} {}_2F_1\left(\frac{1}{6},\frac{1}{6}+\nu;1+\nu;\frac{11-5 \sqrt{5}}{11+5\sqrt{5}}\right)
\end{equation}
\\

For certain values of $\nu$ we get closed form evaluations of (12). An example is $\nu=n+1/2$, $n=0,1,2,\ldots$, which a special case for $n=0$ is
$$
\int^{R^{(-1)}\left(\sqrt[5]{\frac{-11+5\sqrt{5}}{2}}\right)}_{0}f(-q)^4q^{-5/6}R(q)^{5/2}dq=
$$
\begin{equation}
=\frac{27 \Gamma\left(\frac{5}{6}\right) \Gamma\left(\frac{8}{3}\right) \sin\left(\frac{2}{3} \arctan\left(\sqrt{\frac{-11+5 \sqrt{5}}{11+5 \sqrt{5}}}\right)\right)}{2 \sqrt{5 \pi }}
\end{equation} 
We state and prove now the next\\
\\
\textbf{Theorem 2.}\\
If $A=\{1/2,1,1+1/2,2,2+1/2,3,\ldots\}$ and $u(q)=R(q)^{-5}-11-R(q)^5$, then
$$
\int^{u^{(-1)}(x)}_{0}u(q)^{n}\frac{f(-q^5)^5}{f(-q)}dq=
$$
\begin{equation}
=-\int^{x}_{0}\frac{t^{n-1}}{\sqrt{125+22t+t^2}}dt=\textrm{known function when $n\in A$}
\end{equation}
\textbf{Proof.}\\
From [5] we get easily the first equality of (14). For the second equality which is the evaluation part we have ($\rho_1$, $\rho_2$ are the roots of $125+22t+t^2=0$):
$$
\int\frac{t^{n-1}}{\sqrt{125+22t+t^2}}dt
=\int t^{n-1/2}(t-\rho_1)^{-1/2}(t-\rho_2)^{-1/2}dt
$$
is of the form
\begin{equation}
I_{2m}=\int t^{2m}\{(A_1t^2+B_1) (A_2 t^2+B_2)\}^{-1/2}dt
\end{equation}
and reduces to elliptic integrals of the first and second kinds, when $m$ is integer (see [4] and [14]).
\[
\]

Using the same methods of Theorems 1,2 one can prove that
$$
\int^{\rho_1}_{0} x^m(ax^2+bx+c)^n dx= 
$$
\begin{equation}
=\frac{c^n\rho_1^{m+1}\Gamma(m+1)\Gamma(n+1) }{\Gamma(m+n+2)}{}_2F_1\left(m+1,-n;m+n+2;\frac{\rho_1}{\rho_2}\right)
\end{equation}
where $\rho_1<\rho_2$ are the roots of $ax^2+bx+c=0$.    
\[
\]
Another interesting formula is for $0<\nu<1$, (see [5]): 
\begin{equation}
\int^{1}_{0}\left(\frac{f(-q^5)}{f(-q)}\right)^{6\nu+5}f(-q)^4q^{\nu}dq=\frac{-\pi\csc(\nu\pi)}{11^{\nu+1}}{}_2F_{1}\left[\frac{\nu+1}{2},\frac{\nu+2}{2};1;-\frac{4}{121}\right]
\end{equation}
which also can be written as
\begin{equation}
\int^{1}_{0}u(q)^{\nu-5/6}f(-q)^4q^{-5/6}dq=\pi\csc(\nu\pi)11^{\nu-1}{}_2F_{1}\left[\frac{\nu+1}{2},\frac{\nu+2}{2};1;-\frac{4}{121}\right]
\end{equation}
Continuing set $k_r=k$ and $k'_r=k'=\sqrt{1-k^2}$. Then from (7) and the relations  
\begin{equation}
f(-q)=2^{1/3}\pi^{-1/2}q^{-1/24}k^{1/12}k'^{1/3}K(k)^{1/2}
\end{equation}
and
\begin{equation}
\frac{dq}{dk}=\frac{-q\pi^2}{2kk'^2K(k)^2}
\end{equation}
we get 
\begin{equation}
\frac{dR(q)}{dk}=-5^{-1}\cdot2^{1/3}(kk')^{-2/3}R(q)\sqrt[6]{R(q)^{-5}-11-R(q)^5}
\end{equation}
Which is clear that $R(q)$ can be expressed as a function of the singular modulus $k=k_r$.\\
By integration, we get if $q=e^{-\pi\sqrt{r}}$ with $r$ real positive 
\begin{equation}
\frac{3\sqrt[3]{2}}{5}k^{1/3}{}_2F_1\left[\frac{1}{3},\frac{1}{6};\frac{7}{6};k^2\right]=\int^{R(q)}_{0}\frac{dx}{x\sqrt[6]{x^{-5}-11-x^5}}
\end{equation}
\\
Continuing, from (10) replacing were $f$ the Dedekind eta function $\eta$, we get\\
\\
\textbf{Theorem 3.}
\begin{equation}
\pi\int^{+\infty}_{\sqrt{r}}\eta(it/2)^4dt=3\sqrt[3]{2k_r}\cdot{}_2F_1\left[\frac{1}{3},\frac{1}{6};\frac{7}{6};k_r^2\right]=5\int^{R(q)}_{0}\frac{dx}{x\sqrt[6]{x^{-5}-11-x^5}}
\end{equation} 
\[
\]
\textbf{Applications}\\
\textbf{1)} If $r=4$ then from the following evaluation of Ramanujan
$$
R(e^{-2\pi})=-\frac{1+\sqrt{5}}{2}+\sqrt{\frac{5+\sqrt{5}}{2}}
$$
one gets
\begin{equation}
\int^{+\infty}_{2}\eta(it/2)^4dt=5\pi^{-1}\int^{-\frac{1+\sqrt{5}}{2}+\sqrt{\frac{5+\sqrt{5}}{2}}}_{0}\frac{dx}{x\sqrt[6]{x^{-5}-11-x^5}}
\end{equation}
\[
\]
\textbf{2)}
If we set $v(\tau)=R(e^{-\pi\tau})$, $\tau=\sqrt{r}$, then 
$$
-\pi\int^{\tau}_{+\infty}\eta(it/2)^4dt=5\int^{v(\tau)}_{0}\frac{dx}{x\sqrt[6]{x^{-5}-11-x^5}}
$$  
inverting $v_{\tau}$ and derivating we get 
$$
-\pi\eta\left(\frac{iv^{(-1)}(x)}{2}\right)^4\frac{d}{dx}v^{(-1)}(x)=5\frac{1}{x\sqrt[6]{x^{-5}-11-x^5}}
$$
Solving with respect to $\frac{d}{dx}v^{(-1)}(x)$, we get
\begin{equation}
\frac{d}{dx}v^{(-1)}(x)=\frac{-5\pi^{-1}}{\eta\left(i\frac{v^{(-1)}(x)}{2}\right)^4x\sqrt[6]{x^{-5}-11-x^5}}
\end{equation}
and as in Application 1 if we set $x=-\frac{1+\sqrt{5}}{2}+\sqrt{\frac{5+\sqrt{5}}{2}}$, then the next evaluation is valid
$$
\frac{d}{dx}v^{(-1)}\left(-\frac{1+\sqrt{5}}{2}+\sqrt{\frac{5+\sqrt{5}}{2}}\right)=
$$
$$
=-5\pi^{-1}\eta\left(2i\right)^{-4}\sqrt{\frac{1}{10}\left(1+3\sqrt{5}+2\sqrt{10+2\sqrt{5}}\right)}
$$
Where the $\eta\left(2i\right)$ can evaluated from (19) using the singular modulus $k_4=3-2\sqrt{2}$ and the elliptic singular value $K(k_4)$, in terms of algebraic numbers and values of the $\Gamma$ function.
\[
\]
Let $q=e^{-\pi\sqrt{r}}$ and $F(x)$, $m(x)$ be functions defined as\\
\textbf{a)}
$$
x=\int^{F(x)}_{0}\frac{dt}{t\sqrt[6]{t^{-5}-11-t^5}}=\frac{-1}{5}\int^{F(x)^{-5}-11-F(x)^5}_{+\infty}\frac{dt}{t^{1/6}\sqrt{125+22t+t^2}}=
$$
\begin{equation}
=\left(\frac{6}{5}h^{-1/6}F_1\left[\frac{1}{6},\frac{1}{2},\frac{1}{2},\frac{7}{6};-\frac{11-2i}{h},-\frac{11+2i}{h}\right]\right)_{h=F(x)^{-5}-11-F(x)^5}
\end{equation}
and\\
\textbf{b)}
\begin{equation}
x=\pi\int^{+\infty}_{\sqrt{m(x)}}\eta\left(it/2\right)^4dt
\end{equation}
respectively.\\
\\
Using the above functions\\
\\
\textbf{Theorem 4.}\\If $q=e^{-\pi\sqrt{r}}$, $r>0$ then
\begin{equation}
R(q)=F\left(\frac{3}{5}\sqrt[3]{2k_r}\cdot{}_2F_1\left[\frac{1}{3},\frac{1}{6};\frac{7}{6};k_r^2\right]\right).
\end{equation}
\\

Also the equation 
\begin{equation}
\int^{y}_{0}\frac{dt}{t\sqrt[6]{t^{-5}-11-t^5}}=A
\end{equation}
have solution
\begin{equation}
y=R\left(e^{-\pi \sqrt{m(A)}}\right)
\end{equation} 
which also means that the function $F$ have representation  
\begin{equation}
F(x)=R\left(e^{-\pi \sqrt{m(x)}}\right)
\end{equation}
and the solution of the equation $m(x)=r$ is 
\begin{equation}
x=3\sqrt[3]{2k_r}\cdot{}_2F_1\left[\frac{1}{3},\frac{1}{6};\frac{7}{6};k_r^2\right]
\end{equation} 
Also if the differential equation
\begin{equation}
\frac{y'}{y^{\mu}\sqrt{ay^2+by+c}}=\frac{1}{x^{1/6}\sqrt{x^2+22x+125}}
\end{equation}
have solution $y=g(x)$ then $g$ satisfies  
$$
\int^{g(Y(q))}_{0}\frac{dx}{x^{\mu}\sqrt{ax^2+bx+c}}=3\sqrt[3]{2k_r}\cdot{}_2F_1\left[\frac{1}{3},\frac{1}{6};\frac{7}{6};k_r^2\right]
$$
where $Y(q)=R(q^2)^{-5}-11-R(q^2)^{5}$. Hence we get the next\\
\\ 
\textbf{Proposition.}\\
If $g$ satisfies (33), then
\begin{equation}
\int^{g\left(Y\left(e^{-\pi\sqrt{m(x)}}\right)\right)}_{0}\frac{dt}{t^{\mu}\sqrt{at^2+bt+c}}=x
\end{equation}

\section{Some Remarks on a sextic equation}

We define the function $G(x)$ as\\
\begin{equation}
x=\frac{-1}{5}\int^{G(x)}_{+\infty}\frac{dt}{t^{1/6}\sqrt{125+22t+t^2}}
\end{equation}
then $G(x)=F(x)^{-5}-11-F(x)^5$, and the solution of the sextic equation (see [8]): 
\begin{equation}
\frac{b^2}{20a}+bX+aX^2=C_1X^{5/3}
\end{equation}
is
$$
X=X_r=\frac{b}{250a}\left(R(q^2)^{-5}-11-R(q^2)^5\right)=$$
\begin{equation}
=\frac{b}{250a}\cdot G\left(3\cdot5^{-1}\sqrt[3]{2k_{4r}}\cdot{}_{2}F_1\left[\frac{1}{3},\frac{1}{6};\frac{7}{6};k_{4r}^2\right]\right)
\end{equation}
where $j_r=250C_1^3a^{-2}b^{-1}$ is Klein's $j-$invariant. But from 
\begin{equation}
j_r=\frac{16\left(1+14k_{4r}^2+k_{4r}^4\right)^3}{k^{2}_{4r}(1-k_{4r}^2)^4}
\end{equation}
(observe that (38) is always solvable in radicals with respect to $k_{4r})$, we lead to the following\\
\\
\textbf{Theorem 5.}\\
A solution of (36) is
\begin{equation}
X=\frac{b}{250a}G\left(3\sqrt[3]{2}t^{1/6}{}_2F_{1}\left[\frac{1}{3},\frac{1}{6};\frac{7}{6};t\right]\right)
\end{equation}
where the $t$ is given from 
\begin{equation}
250C_1^3a^{-2}b^{-1}=\frac{16(1+14t+t^2)^3}{t(1-t)^4}.
\end{equation}
Relation (39) is always solvable with respect to $t$. If we make the change of variable 
$$
t=\left(\frac{1-\sqrt{1-x^2}}{1+\sqrt{1-x^2}}\right)^2 , 
$$ 
then (40) becomes solvable in radicals with respect to $x$ and hence to $t$ also.\\
\\

From Theorem 5 and the definition of $G(x)$ it is clear that
$$ 
-\frac{1}{5}\int^{250X_rab^{-1}}_{+\infty}\frac{dx}{x^{1/6}\sqrt{125+22x+x^2}}=5^{-1}4^{-1/3}B\left(b_r,\frac{1}{6},\frac{2}{3}\right)
$$
where 
$$
250C_1^3a^{-2}b^{-1}=\frac{16(1+14b_r+b_r^2)^3}{b_r(1-b_r)^4} . 
$$
Setting $a=1$, $b=250$ we get the next\\
\\
\textbf{Proposition 1.}\\
We define the function $\beta^*=\beta^*_r$ to be root of the equation 
\begin{equation}
\sqrt{\frac{B\left(1-x,\frac{1}{6},\frac{2}{3}\right)}{B\left(x,\frac{1}{6},\frac{2}{3}\right)}}=\sqrt{r}
\end{equation}
Then $x=\beta^*_r$ is algebraic when $r\in Q^{*}_{+}$ and further for any $b_r$
\begin{equation}
B\left(b_r,\frac{1}{6},\frac{2}{3}\right)=\sqrt[3]{4}\int^{+\infty}_{X_r}\frac{dx}{x^{1/6}\sqrt{125+22x+x^2}}  
\end{equation}
with
\begin{equation}
C_1^3=\frac{16\left(1+14b_r+b_r^2\right)^3}{b_r\left(1-b_r\right)^4}  
\end{equation}
and $X_r$ is solution of
\begin{equation}
3125+250X_{r}+X_{r}^2=C_1X_{r}^{5/3} .
\end{equation}
Also if $r_1=4^{-1}k^{(-1)}\left(\beta_r^{1/2}\right)$ then $\beta_r=b_{r_1}=k_{4r_1}^2$ and $C_1=j_{r_1}^{1/3}$ and 
\begin{equation}
X_{r_1}=R\left(e^{-2\pi\sqrt{r_1}}\right)^{-5}-11-R\left(e^{-2\pi\sqrt{r_1}}\right)^5
\end{equation}
\\
\textbf{Note.} From the above Proposition we get also that for any $X>0$ we can find $b=\theta(X)$, where $\theta(X)$ is an algebraic function of $X$ and can be calculated explicit from
\begin{equation}
\int^{+\infty}_{X}\frac{dx}{x^{1/6}\sqrt{125+22x+x^2}}=\frac{1}{\sqrt[3]{4}}B\left(b,\frac{1}{6},\frac{2}{3}\right)
\end{equation} 
\\
\textbf{Theorem 6.}\\
If $F$ is arbitrary function and the algebraic singular equation
\begin{equation}
\frac{F(1-x)}{F(x)}=\sqrt{r}
\end{equation}
have root $x=\alpha_r$ then we can define $r_0$ as
\begin{equation}
r_0=\frac{K\left(\sqrt{1-\alpha_r}\right)^2}{K\left(\sqrt{\alpha_r}\right)^2} .
\end{equation} 
If the computated $r_0$ gives the value of $j_{r_0}$ in radicals then we can evaluate the value of $\alpha_r$ from the solvable equation
\begin{equation}
\frac{256(x+(1-x)^2)^3}{x^2(1-x)^2}=j_{r_0} .
\end{equation}
More precisely one algebraic solution will be $x=\alpha_r$.\\ 
\\

The above elementary theorem is for numerical purposes and is a complete change of base from $F$ to $K$ and the evaluation of $\alpha_r$ by the class invariant $j_{r_0}$.\\ 
Note also that the theorem holds for every $F$ not necessary modular base.\\
\\
\textbf{Theorem 7.}\\
If $$
X(r)=R\left(e^{-2\pi\sqrt{r}}\right)^{-5}-11-R\left(e^{-2\pi\sqrt{r}}\right)^{5}
$$
then 
$$
X'(r)=-\pi\frac{\eta(i\sqrt{r})^4}{\sqrt{r}}X(r)^{1/6}\sqrt{125+22 X(r)+X(r)^2}
$$
\textbf{Proof.}\\
From Theorem 3 and definition (b) relation (27), we have that 
$$
m(r)=-2^{5/3}B\left(k_r^2,1/6,2/3\right)
$$
and
$$
\frac{d}{dr}B\left(k_r^2,1/6,2/3\right)=\pi\sqrt[3]{4}\frac{\eta\left(i\sqrt{r}/2\right)^4}{2\sqrt{r}}
$$
But 
$$
\frac{-1}{5}\int^{X(r/4)}_{+\infty}\frac{dt}{t^{1/6}\sqrt{125+22t+t^2}}=\frac{1}{5\sqrt[3]{4}}B\left(k_r^2,1/6,2/3\right)
$$
Differentiating we get the result.\\
\\

\textbf{3)} Continuing the investigation of equation (36) we have if $X=X_r=Y(q)$, then  
$$
C_1Y(q)^{5/3}=aY(q)^2+bY(q)+b^2/(20a)$$ hence, from (9) with 
$$
j_r=250C_1^3a^{-2}b^{-1}
$$ 
and 
$$
Y(q)=R(q^2)^{-5}-11-R(q^2)^5
$$
then
$$
\frac{d}{dx}\left(\frac{3x^{8/3}}{25000}F_1\left[\frac{8}{3},1,1,\frac{11}{3};\frac{-x}{25(5-2\sqrt{5})},\frac{-x}{25(5+2\sqrt{5})}\right]\right)_{x=Y(q)}=j_r^{-1/3} 
$$
Hence if $j_r=j(q)$, then
\begin{equation}
\int^{x}_0 j\left(Y^{(-1)}(w)\right)^{1/3}dw=\frac{3x^{8/3}}{25000}F_1\left[\frac{8}{3},1,1,\frac{11}{3};\frac{-x}{25(5-2\sqrt{5})},\frac{-x}{25(5+2\sqrt{5})}\right] 
\end{equation}
But as we will see in Section 3 bellow, we have $\nu=1$ and $P_n(\nu,x)=\sum^{n}_{k=0}x^k$, for all numbers $n$. 
Hence if 
$$
\phi(x)={}_2F_{1}\left[1,\frac{8}{3};\frac{11}{3};x\right]
$$
then
$$
\int^{x}_{0}j\left(Y^{(-1)}(w)\right)^{1/3}dw=
$$
\begin{equation}
=\frac{3x^{8/3}}{25000(x-1)}\left[-\phi\left(\frac{5+2\sqrt{5}}{125}x\right)+x\phi\left(\frac{5+2\sqrt{5}}{125}x^2\right)\right] .
\end{equation}
The function $\phi$ can evaluated using elementary functions. For $x=1$ we get the next special value
$$
\int^{1}_{0}j\left(Y^{(-1)}(w)\right)^{1/3}dw=\frac{100}{719} \left(60+\sqrt{5}\right)-5\cdot {}_2F_1\left[1,\frac{8}{3};\frac{11}{3};\frac{1}{125} \left(5+2 \sqrt{5}\right)\right]
$$

\section{Modular equations and singular values of Beta functions}

From relations (8) and (9) we have
$$
\int \frac{x^{\mu}}{(ax^2+bx+c)^{\nu}}dx=c^{-\nu}\frac{x^{\mu+1}}{\mu+1}\sum^{\infty}_{n=0}\frac{(\mu+1)_n}{(\mu+2)_n}\left(\sum^{n}_{l=0}\frac{(\nu)_l(\nu)_{n-l}}{l!(n-l)!}\frac{1}{\rho_2^l\rho_1^{n-l}}\right)x^n
$$ 
If we denote $P_n(x)$ as
\begin{equation}
P_n(\nu,x)={}_2F_1\left[-n,\nu;1-n-\nu;x\right]
\end{equation}
then
$$
\sum^{n}_{l=0}\frac{(\nu)_l(\nu)_{n-l}}{l!(n-l)!}\frac{1}{\rho_2^l\rho_1^{n-l}}=\rho_1^{-n}\left(\nu\right)_n{}_2F_1\left[-n,\nu;1-n-\nu;\frac{\rho_1}{\rho_2}\right]=
$$
$$
=\rho_1^{-n}\left(\nu\right)_nP_{n}\left(\nu,\frac{\rho_1}{\rho_2}\right)
$$   
hence if $C_{n,k}=\textrm{Binomial}[n,k]=\left(^n_k\right)$, and
$$
\int^{x}_{0} \frac{t^{\mu}}{(at^2+bt+c)^{\nu}}dt=c^{-\nu}x^{\mu+1}\sum^{\infty}_{n=0}\frac{\left(\nu\right)_nP_n\left(\nu,\frac{\rho_1}{\rho_2}\right)}{(n+\mu+1)\rho^n_1}\frac{x^n}{n!}=
$$
\begin{equation}
=c^{-\nu}x^{\mu+1}\sum^{\infty}_{n=0}\frac{\left(\nu\right)_n}{(n+\mu+1)\rho^n_1}\left[\sum^{n}_{\lambda=0}\frac{C_{n,\lambda}(\nu)_{\lambda}}{(1-n-\nu)_{\lambda}}\left(\frac{-\rho_1}{\rho_2}\right)^{\lambda}\right]\frac{x^n}{n!}
\end{equation}
From the above we obtain the next Ramanujan-type (like $\pi$) formula\\
\\
\textbf{Theorem 8.}\\
Let $\rho_1$, $\rho_2$ be the roots of $ax^2+bx+c=0$, then
$$
\int^{x}_{0} \frac{t^{\mu}}{(at^2+bt+c)^{\nu}}dt=
$$
\begin{equation}
=c^{-\nu}x^{\mu+1}\sum^{\infty}_{n=0}\left[\sum^{n}_{\lambda=0}\frac{C_{n,\lambda}\cdot(\nu)_{\lambda}}{(1-n-\nu)_{\lambda}}\left(\frac{-\rho_1}{\rho_2}\right)^{\lambda}\right]\frac{\left(\nu\right)_n}{n!}\frac{(x/\rho_1)^n}{n+\mu+1}
\end{equation}
\\
\textbf{Examples.}\\
\textbf{1)} For $m=0$, $\nu=1/2$, and $a=b=c=1$, we get
\begin{equation}
\log\left(1+\frac{2}{\sqrt{3}}\right)=\sum^{\infty}_{n=0}\left[\sum^{n}_{\lambda=0}\frac{C_{n,\lambda}\cdot\left(\frac{1}{2}\right)_{\lambda}}{(\frac{1}{2}-n)_{\lambda}}\left(\frac{1-i\sqrt{3}}{2}\right)^{\lambda}\right]\frac{\left(\frac{1}{2}\right)_n}{(n+1)!}\left(\frac{2i}{\sqrt{3}-i}\right)^n
\end{equation}
\textbf{2)} Let now $\mu=-1/2$, $\nu=1/2$, $a=1$ and $\rho_1=-\rho_2=p$, $x=p^2$. Then
$$
\frac{1}{p}\int^{p^2}_{0}\frac{dt}{\sqrt{t(1-t/p)(1+t/p)}}=2\cdot{}_2F_{1}\left[\frac{1}{2},\frac{1}{4};\frac{5}{4},p^2\right]=\frac{\sqrt{p}}{2}B\left(p^2,\frac{1}{4},\frac{1}{2}\right)=
$$
$$
=\sum^{\infty}_{n=0}\left(\sum^{n}_{\lambda=0}(-1)^{\lambda}\frac{C_{n,\lambda}\left(\frac{1}{2}\right)_{\lambda}}{\left(\frac{1}{2}-n\right)_{\lambda}}\right)\frac{\left(\frac{1}{2}\right)_n}{p^n\left(n+\frac{1}{2}\right)n!}
$$
where 
\begin{equation}
B(x,\alpha,\beta)=\int^{x}_{0}t^{\alpha-1}(1-t)^{\beta-1}dt
\end{equation}
\\
\textbf{Conjecture.}\\
The functions $\sqrt{B(z,\alpha,\beta)}$, where\\ $\{\alpha,\beta\}\in\{\{\frac{1}{2},\frac{1}{2}\},\{\frac{1}{2},\frac{1}{3}\},\{\frac{1}{2},\frac{1}{4}\},\{\frac{1}{2},\frac{1}{6}\},\{\frac{1}{3},\frac{1}{3}\},\{\frac{1}{3},\frac{1}{6}\},\{\frac{1}{6},\frac{1}{6}\}\}
$ 
are modular bases.\\
\\
These are modular bases which I manage to find with a quick view. I suppose that exist more of them. Also these bases are represented from Gauss hypergeometric functions. I think that these bases are not like the cubic the 4th the fifth modular bases which related with the complete elliptic integral of the first kind $K$ (see [6]).\\
\\   
\textbf{Example 3.}\\
Consider the equation
$$
\frac{\sqrt{B\left(1-\beta_r,\frac{1}{6},\frac{1}{6}\right)}}{\sqrt{B\left(\beta_r,\frac{1}{6},\frac{1}{6}\right)}}=\sqrt{r}
$$
then some solutions extracted with numerical methods of Mathematica are 
$$
\beta_2=\frac{1}{4}\left(2-\sqrt{3}\right)
$$
$$
\beta_3=\frac{1}{4}\left(2-\sqrt{3 \left(3-\sqrt{3}\right)}\right)
$$
$$
\beta_{3/2}=\frac{1}{8} \left(4-\sqrt{-9+9 \sqrt{5}-3 \sqrt{150-66 \sqrt{5}}}\right)
$$
$$
\beta_4=\frac{1}{8} \left(4-\sqrt{-9+9 \sqrt{5}+3 \sqrt{150-66 \sqrt{5}}}\right)
$$
$$
\beta_5=\frac{1+3\cdot 2^{1/3}-3\cdot 2^{2/3}}{8+4\cdot \sqrt{3\cdot \left(1-2^{1/3}+2^{2/3}\right)}}
$$
and have the property 
\begin{equation}
\frac{B\left(\beta_{n^2r},\frac{1}{6},\frac{1}{6}\right)}{B\left(\beta_{r},\frac{1}{6},\frac{1}{6}\right)}=\textrm{rational} ,
\end{equation}
when $n,r\in\bf N\rm$.\\
Consider the function\\
\\
\textbf{Definition.} $$B_\alpha(x):=\sqrt{B(x,\alpha,\alpha)}=\sqrt{\int^{x}_{0}(t-t^2)^{\alpha-1}dt},\textrm{ where }0<a<1\textrm{ and }x>0$$ 
\\

Then by changing of variable $t\rightarrow 1-w$, we have  
$$
B^2_{\alpha}(x)+B^2_{\alpha}(1-x)=\int^{1}_{0}(w(1-w))^{\alpha-1}dw=\frac{\Gamma(\alpha)^2}{\Gamma(2\alpha)}
$$
Setting now $x=\beta_r$, where 
\begin{equation}
\frac{B_{\alpha}(1-\beta_r)}{B_{\alpha}(\beta_r)}=\sqrt{r}
\end{equation}
we conclude that 
\begin{equation}
B_{\alpha}(\beta_r)=\sqrt{\frac{\Gamma(\alpha)^2}{\Gamma(2\alpha)(r+1)}}
\end{equation}
From Theorem 8 with $a=1$, $b=-2$, $c=1$, $\mu=-5/6$, $\nu=5/12$ and $r=5$, $x=\beta_5$, we have the next formula for the constant $B\left(\beta_5,\frac{1}{6},\frac{1}{6}\right)=\frac{\Gamma\left(\frac{1}{6}\right)}{\sqrt{6\Gamma\left(\frac{1}{3}\right)}}$ :\\
\begin{equation}
\frac{\Gamma\left(\frac{1}{6}\right)^2}{6\Gamma\left(\frac{1}{3}\right)}
=\sum^{\infty}_{n=0}\left(\sum^{n}_{\lambda=0}(-1)^{\lambda}\frac{C_{n,\lambda}\left(\frac{5}{12}\right)_{\lambda}}{\left(\frac{7}{12}-n\right)_{\lambda}}\right)\frac{\left(\frac{5}{12}\right)_n\left(\beta_{5}\right)^{n+1/6}}{n!\left(n+\frac{1}{6}\right)}
\end{equation}
\\
In general holds the following very interesting theorem\\
\\
\textbf{Theorem 9.}\\
If $\beta_r$ is solution of (58) then 
\begin{equation}
\frac{\Gamma\left(\alpha\right)^2}{(r+1)\Gamma\left(2\alpha\right)}=\sum^{\infty}_{n=0}\left(\sum^n_{\lambda=0}(-1)^{\lambda}\frac{C_{n,\lambda}\left(\frac{1-\alpha}{2}\right)_{\lambda}}{\left(\frac{1+\alpha}{2}-n\right)_{\lambda}}\right)\frac{\left(\frac{1-\alpha}{2}\right)_{n}\beta_r^{n+\alpha}}{n!(n+\alpha)}
\end{equation}
\\
\textbf{Example 4.}\\
Let 
$$
\psi(x):=\sqrt{2\arcsin\left(x\right)}
$$
Then the equation
\begin{equation}
\frac{\psi(1-s_r)}{\psi(s_r)}=\sqrt{r}\textrm{ , }r \in \textbf{R}^{*}_{+} 
\end{equation}
is equivalent to
\begin{equation}
2i+y^{-1}-y+y^{-r}-y^r=0
\end{equation}
with $y=is+\sqrt{1-s^2}$, $s=s_r$.\\
\textbf{Proof.}\\
Equation (62) is written as
$$
\arcsin(1-x)=r\arcsin(x)\textrm{ or }\sin(r^{-1}\arcsin(1-x))=x
$$
which is similar with Tchebychev polynomials $T_r(x)=\cos(r\arccos(x))$, $r>0$.\\Using the identity 
$$
(\cos(x)+i\sin(x))^r=\sum^{\infty}_{n=0}C_{r,n}i^n\sin(x)^n\left(\sqrt{1-\sin(x)^2}\right)^{r-n}=
$$
$$
=\cos(nx)+i\sin(nx)
$$
and setting $\sin(x)=w$ one can arrive to
\begin{equation}
\left(\sqrt{1-w^2}-iw\right)^r-\left(\sqrt{1-w^2}+iw\right)^r+2i(1-w)=0
\end{equation}
which is (63).\\
\\

Set now $\xi=\frac{1}{2} \left(-i-\sqrt{3}\right)$ and for every $r>0$ 
$$
x_r=2^{-1/r} \left(\frac{1-2 i \xi-\xi^2+\sqrt{1-4 i \xi-2 \xi^2+4 i \xi^3+\xi^4}}{\xi}\right)^{\frac{1}{r}}
$$
If
$$
x=\frac{-i(1-x_r^2)}{2x_r},
$$
and solve the equation $x=\sin(t)$ with respect to $t$ (i.e $t=\arcsin(x)$). Finally we get that $t=\sin\left(\pi/(4+2r)\right)$ is solution of
$$
\frac{\psi(1-2t^2)}{\psi(t)}=\sqrt{r}
$$

Another interesting note is for the function 
\begin{equation}
m(R)=-\frac{1}{4}e^{\frac{-i\pi}{1+R}}\left(-1+e^{\frac{i\pi}{1+R}}\right)^2=\sin\left(\frac{\pi}{2(R+1)}\right)^2
\end{equation} 
hold the following\\
\\
\textbf{Example 5.}\\
Let 
$$
\psi^{*}(x)=\sqrt{B_{1/2}\left(x,\frac{1}{2},\frac{1}{2}\right)}=\sqrt{\arcsin\left(x^{1/2}\right)}
$$ 
then $m(R)$ is the solution of
\begin{equation}
\frac{\psi^{*}(1-x)}{\psi^{*}(x)}=\sqrt{R}\textrm{ , }R \in \textbf{R}^{*}_{+} 
\end{equation}
and holds the following modular equation
\begin{equation}
m(R+1)=\frac{1-\sqrt{1-m\left(\frac{R}{2}\right)}}{2}
\end{equation}
\\  
\textbf{Example 6.}\\
If $m(R)$ is given form (67), then
\begin{equation}
\frac{\pi}{R+1}=\sum^{\infty}_{n=0}\frac{\left(\frac{1}{4}\right)_n\left(\frac{3}{4}\right)_{-n}}{\left(\frac{1}{2}\right)_{-n}n!}\frac{m(R)^{n+1/2}}{n+1/2}
\end{equation}
\textbf{Proof.}
\begin{equation}
B\left(m(R),\frac{1}{2},\frac{1}{2}\right)=\frac{\pi}{R+1}=\sum^{\infty}_{n=0}\left(\sum^{n}_{\lambda=0}\frac{(-1)^{\lambda} C_{n,\lambda}\left(\frac{1}{4}\right)_{\lambda}}{\left(\frac{3}{4}-n\right)_{\lambda}}\right)\frac{\left(\frac{1}{4}\right)_n m(R)^{n+1/2}}{n!\left(n+\frac{1}{2}\right)}
\end{equation}
The above is a Ramanujan-type $\pi$ formula of arbitrary large precision since any higher order values of $m(R)$ can always evaluated from (65) in radicals. Moreover is 
$$
\sum^{\infty}_{n=0}\frac{\left(\frac{1}{4}\right)_n\left(\frac{3}{4}\right)_{-n}}{\left(\frac{1}{2}\right)_{-n}n!}\frac{m(R)^{n+1/2}}{n+1/2}=2\arcsin\left(\sqrt{m(R)}\right)
$$ \\
\\
\textbf{Note.} Formulas similar to (67) can given if we consider the expansion of
$$
\sin(nt)=P_n\left(\sqrt{1-\sin(t)^2},\sin(t)\right)
$$
For example with $n=6$ we get 
$$
p_{6}(y)=P_6(x,y)=6x^5y-20x^3y^3+6xy^5\textrm{ , }x=\sqrt{1-y^2}$$ and for this $n$ the equation  $p_6(t)=P_6(\sqrt{1-t^2},t)=a$ is solvable. Hence one can get a formula 
$$
t=\sin\left(\frac{\theta}{6}\right)=p_6^{(-1)}\left(\sin(\theta)\right)
$$

\[
\]

\centerline{\bf References}\vskip .2in

\noindent

[1]: M. Abramowitz and I.A. Stegun. 'Handbook of Mathematical Functions'. Dover Publications, New York.(1972)

[2]: G.E. Andrews. 'Number Theory'. Dover Publications, New York.(1994)
 
[3]: G.E. Andrews, Amer. Math. Monthly, 86, 89-108.(1979)

[4]: J.V. Armitage, W.F. Eberlein. 'Elliptic Functions'. Cambridge University Press.(2006)

[5]: N.D. Bagis and M.L. Glasser. 'Integrals related with Rogers Ramanujan continued fraction and q-products'. arXiv:0904.1641v1 [math.NT].(2009)

[6]: N.D. Bagis and M.L. Glasser. 'Conjectures on the evaluation of alternative modular bases and formulas approximating 1/$\pi$'. Journal of Number Theory. (Elsevier).(2012)

[8]: N.D. Bagis. 'On a General Sextic Equation Solved by the Rogers Ramanujan Continued Fraction'. arXiv:1111.6023v2[math.GM].(2012) 

[9]: B.C. Berndt. 'Ramanujan`s Notebooks Part I'. Springer Verlag, New York.(1985)

[10]: B.C. Berndt. 'Ramanujan`s Notebooks Part II'. Springer Verlag, New York.(1989)

[11]: B.C. Berndt. 'Ramanujan`s Notebooks Part III'. Springer Verlag, New York.(1991)

[12]: I.S. Gradshteyn and I.M. Ryzhik. 'Table of Integrals, Series and Products'. Academic Press.(1980)  

[13]: N.N. Lebedev. 'Special Functions and their Applications'. Dover Pub. New York.(1972)

[14]: E.T. Whittaker and G.N. Watson. 'A course on Modern Analysis'. Cambridge U.P.(1927)

\end{document}